\documentclass{article}
\usepackage{amssymb,amsmath,amsfonts,cite,bm}
\usepackage{graphicx}
\usepackage{cite}

\newtheorem{theorem}{Theorem}[section]
\newtheorem{lemma}[theorem]{Lemma}

\title
{Realization of Arbitrary Hysteresis by a Low-dimensional Gradient Flow}

\begin{document}

\date{}

\maketitle

\centerline{\scshape Dmitrii Rachinskii\footnote{The author was supported by National Science Foundation grant DMS-1413223.}}
\medskip
{\footnotesize
\centerline{Department of Mathematical Sciences}
\centerline{The University of Texas at Dallas, Richardson, TX, USA}
}
\bigskip




\bigskip

\begin{abstract}
We consider gradient systems with an increasing potential that depends on a scalar parameter. As the parameter is varied, critical points of the potential can be eliminated or created through saddle-node bifurcations causing the system to transit from one stable equilibrium located at a (local) minimum point of the potential to another minimum along the heteroclinic connections. These transitions can be represented by a graph. We show that any admissible graph has a realization in the class of two dimensional gradient flows. The relevance of this result is discussed in the context of genesis of hysteresis phenomena. The Preisach hysteresis model is considered as an example.
\end{abstract}

\section{Introduction}
Hysteresis is a complex type of relationship between variables. The term was first introduced by James Ewing in his studies of properties of magnetic materials and may be translated as ``to lag behind'' from ancient Greek. The key feature of hysteresis is the dependence of the value of an output variable on some past values of the input variable, that is a type of memory. This dependence is mediated by the internal state that changes in response to variations of the input and is affected by the history of these variations (hence, to predict future outputs, one needs to know either the current internal state or the history of input variations). The term hysteresis is associated with a specific type of memory
that can be characterized as persistent memory (as opposed to fading or scale-dependent memory such as described by the convolution integrals, delayed arguments, etc.). A mathematical idealization that is used to describe hysteresis is a so-called rate-independent input-output operator, that is an operator that commutes with the group of all increasing transformations of time \cite{Vis}. A consequence of the rate-independence (in case of scalar-valued inputs) is that only certain (local) maximum and minimum values of the input achieved in the past can affect the future. 

\begin{figure}
    \centering
        \includegraphics[width=0.33\textwidth]{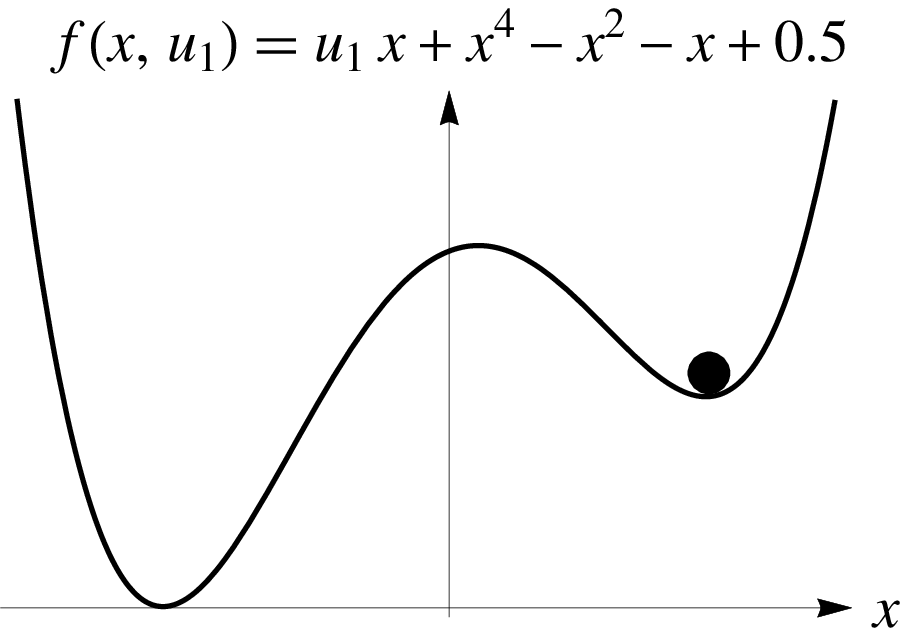} \qquad\quad \includegraphics[width=0.33\textwidth]{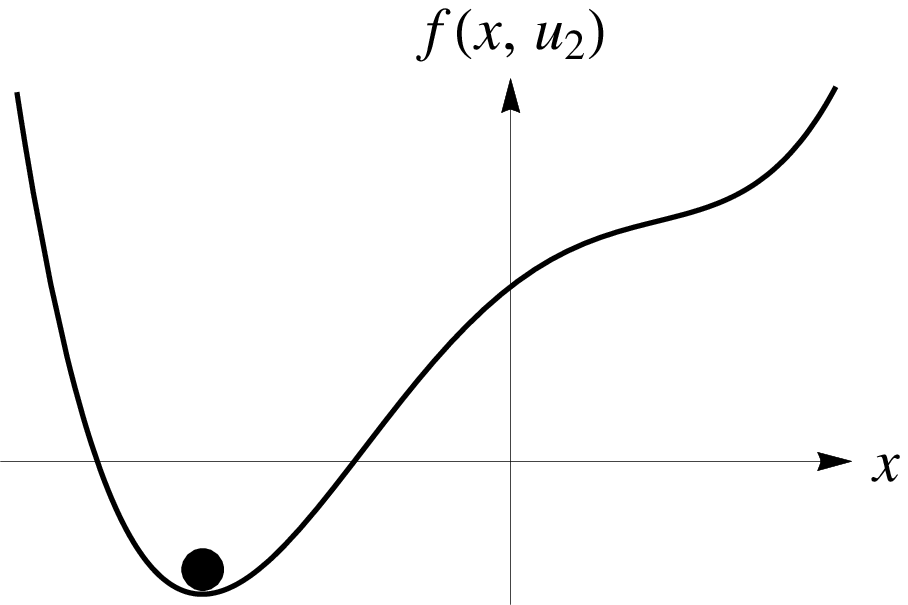}
    \caption{Double well potential depending on a control parameter $u$. Black point indicates the state of the system. When $u$ increases, the right minimum collides with the local maximum of the potential
    and disappears. As a result of this bifurcation, the black point that was sitting in the right minimum (left panel) makes a transition to the remaining minimum (right panel). } \label{Fig1}
\end{figure}

\begin{figure}
    \centering
         \begin{minipage}{0.45\textwidth}\includegraphics[width=1\textwidth]{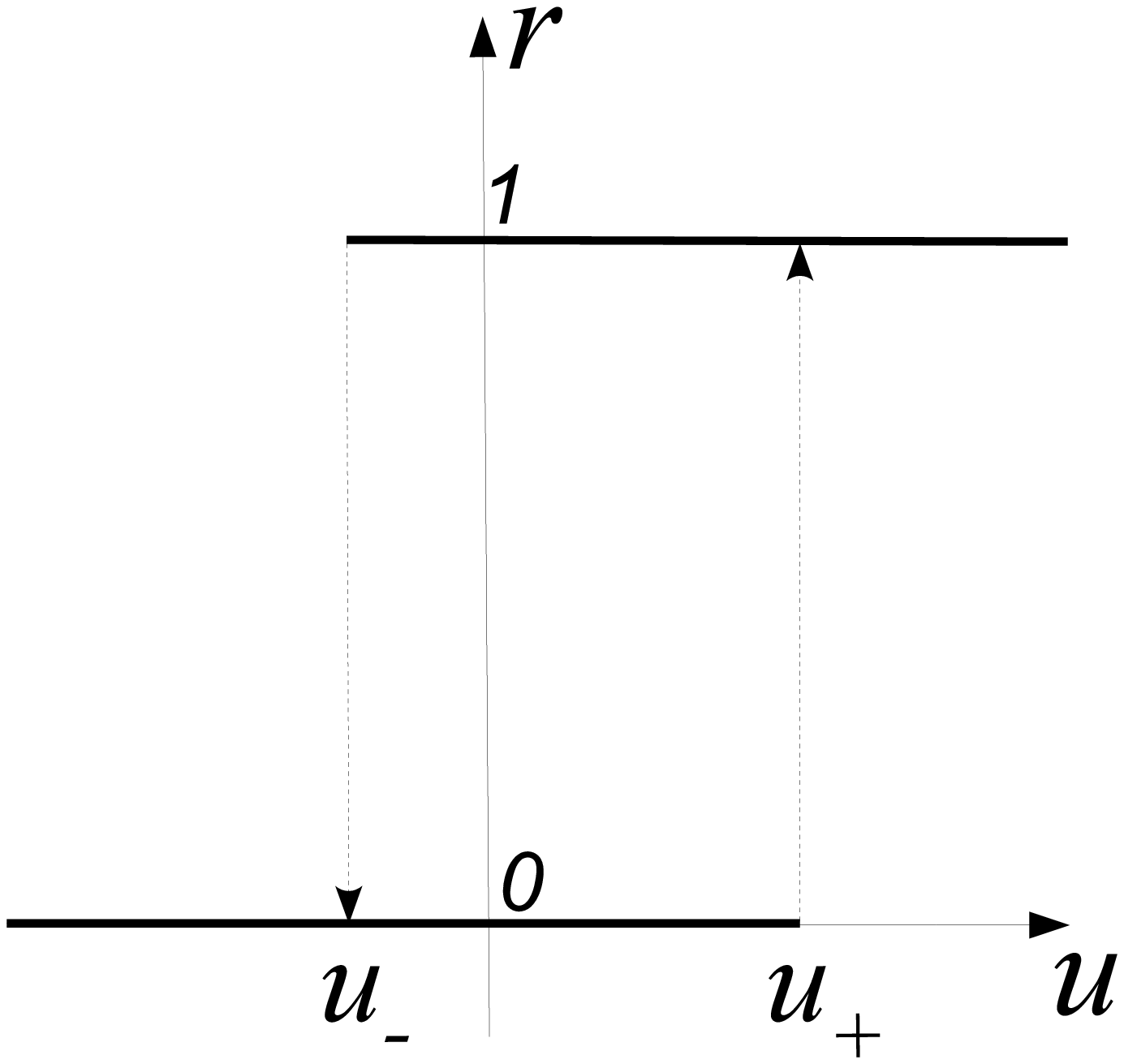}\end{minipage} \qquad\quad  \begin{minipage}{0.3\textwidth}\includegraphics[width=1\textwidth]{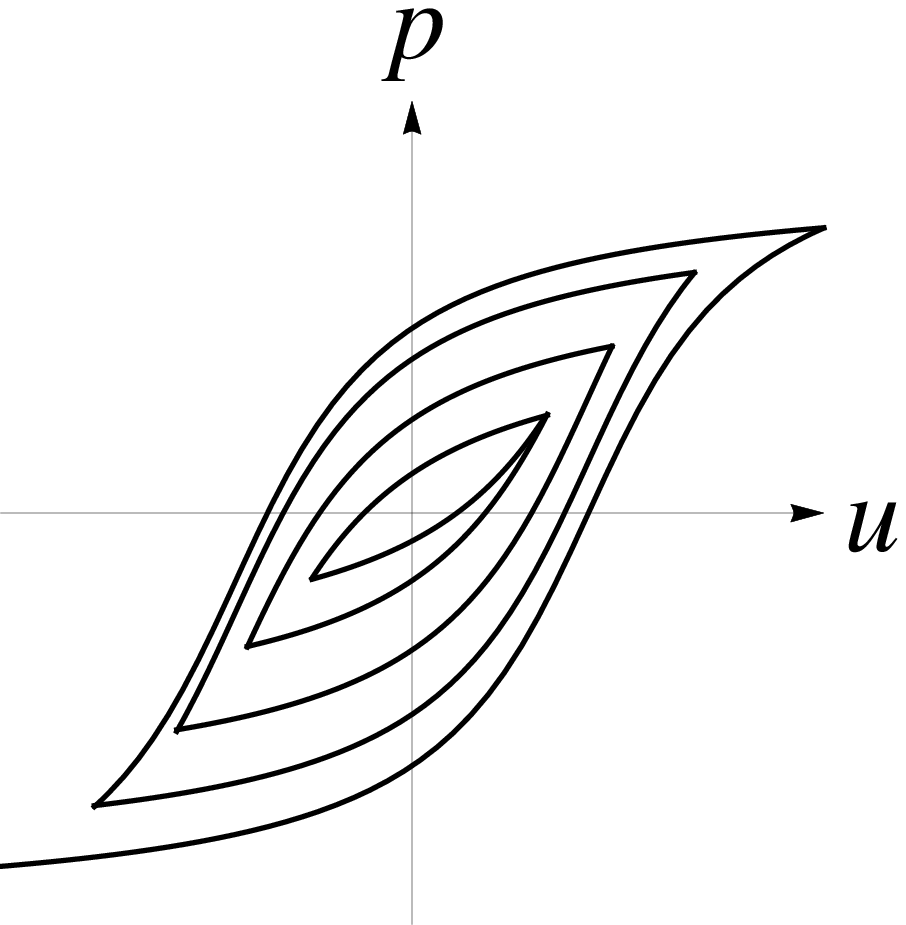}\end{minipage}
    \caption{{\bf Left:} Two-state non-ideal relay switches from state $r=0$ to state $r=1$ when the input $u$ exceeds the threshold value $u_+$
    and from state $r=1$ to state $r=0$ when $u$ goes below the threshold value $u_-$ that satisfies $u_-<u_+$. {\bf Right:} Hysteresis curves typical of complex hysteresis relationship between
    input $u$ and output $p$ such as in the Preisach model. The picture corresponds to an oscillating input with decreasing amplitude that becomes periodic after a few oscillations. } \label{Fig2}
\end{figure}

Hysteresis occurs in many phenomena, for example, it is encountered in magnetism \cite{May, Mal}, plasticity \cite{Lub,mroz}, friction \cite{Lam, Rud}, mechanical damage and fatigue
\cite{Mic}, constitutive relationships of ``smart'' materials (shape memory alloys, piezo-electric and magnetostrictive materials \cite{Kuh, Cir}), porous media filtration \cite{Par,Rah,Bot}, phase transitions \cite{BroSpr,Dah} and engineering (thermostat, non-ideal switch and backlash nonlinearity are usual examples \cite{Mih,Gur}).
The simplest example of hysteresis is probably a bi-stabile system (in dynamical systems theory the two terms, bi-stability and hysteresis, are sometimes used almost as synonyms).
Suppose that at an equilibrium the system minimizes some potential function (such as Helmholtz free energy density in Landau's theory of phase transitions and critical phenomena \cite{Lan}). Fig.~\ref{Fig1} shows a double well potential perturbed by a linear term. Let us consider the coefficient $u$ of this term as a control parameter (input) of the system. Suppose the system occupies the state corresponding to the right minimum point of the potential at some time. As $u$ varies, the profile of the potential changes, and the right minimum can be eliminated, in which case the system has to make a transition to the other (remaining) minimum point, see Fig.~\ref{Fig1}. A similar transition from the left to the right minimum point occurs at a different threshold value of the control parameter $u$, at which the left minimum point gets eliminated. Thus, we observe co-existence of two equilibrium states for some interval $[u_-,u_+]$ of the input values with transitions from one state to the other happening when $u$ reaches the value $u_-$ and the opposite transitions happening at the other end $u_+$ of the bi-stability interval. This is exactly a description of the simple hysteretic switch known as the non-ideal two-threshold relay, see Fig.~\ref{Fig2} (left).

It is important to note that the applicability of the bi-stable relay model hinges on certain assumptions about characteristic time scales of several processes. In statistical physics, the global minimum of the energy potential corresponds to a persistent stable equilibrium, while local (relative) minima are interpreted as metastable states that can persist for limited characteristic time only before decaying eventually to the global minimum due to thermal fluctuations. The non-ideal relay model assumes that the characteristic time that passes between the transitions of the system from one state to another due to variations of the input $u$ is much shorter than the characteristic time that the system would spend in a metastable state before thermal fluctuations would guide it to the absolute minimum. On the other hand, the process of relaxation to a new state after the system has left some metastable state (that has been eliminated by a variation of the input) is assumed to be much faster than the input process, that is transitions are almost instantaneous compared to the time that the system spends in the stable and metastable states. Let us also remark that a mechanisms of switching between states of a bi-stable system is not necessarily related to minimization of some functional. Two further important examples of multi-stable systems are presented by slow-fast systems \cite{Oma} and systems with feedback (such as, for instance, in Weiss mean field model of magnetization \cite{Wei}).

Fig.~\ref{Fig1} presents the simple hysteresis effect. More complex hysteresis phenomena (that have been well known in magnetism, plasticity and sorption for decades) are described by input-output diagrams with a continual family of ascending and descending curves and nested hysteresis loops, see Fig.~\ref{Fig2} (right). A schematic input-output diagram on the right panel of Fig.~\ref{Fig2} may represent a complex dependence of magnetization of a ferromagnetic material on the applied magnetic field, or a stress-strain relationship in an elasto-plastic body, or moisture content of a porous medium vs applied pressure, etc. Various models of complex hysteresis relationships have been proposed and used by physicists and engineers.  A few most prominent models include the Preisach model of magnetic materials, the Prandtl-Ishlinskii model of plasticity and friction and the Ising spin-interaction model of phase transitions in statistical physics (that have been also adapted to model sorption hysteresis \cite{OKa,Hyd,Iye,energy}; damage accumulation \cite{Ryc}; constitutive laws of ``smart'' materials; hysteresis in economics \cite{Fin,Pok, Cro}, social dynamics and population biology \cite{Mel,Bio,Bio1,Bio2}; see \cite{Sci, Siam} for further examples). Mathematical tools that have been developed in the area of modeling hysteresis phenomena and systems with hysteretic components are diverse and include the method of hysteresis operators \cite{KraPok}, differential inclusions \cite{Kre}, variational inequalities \cite{KreBoo}, variational approach based on $\Gamma$-convergence \cite{Mie} and switched systems \cite{Ast}.

The above mentioned models of hysteresis are phenomenological and either describe hysteresis effects on a macroscopic level (for example, the Preisach and Prandtl-Ishlinskii models) or use an extremely simplified microscopic model of reality (such as in the Ising model). Derivation of macroscopic models of hysteresis from first principles, or microscopic physics, remains a 
daunting challenge \cite{Sci}. One idea is that in a complex energy landscape with a large number of metastable states, transitions between these states caused by variation of the control parameter (such as shown in Fig.~\ref{Fig1} for a simple double well potential) may result in a complex hysteresis relationship between macroscopic variables.  Indeed, this idea has been entertained by researchers in hysteretic systems for a period of time. However, to the best of our knowledge, no precise results validating this idea have been obtained so far. A question that may be asked in this relation is whether any given ``map'' of transitions between metastables states can be realized by a gradient system that transits from one equilibrium to another in response two variations of the input parameter.

In this paper, we give a positive answer to this question under the assumptions that (a) the hysteretic system has a discrete set of states (and hence can be described as a directed graph); and, (b) the input of the system is a scalar variable. It turns out that any system of this type can be realized by a gradient system with two degrees of freedom. We do not impose restrictions on the functional of the gradient system. Whether it is possible to use functionals from a certain class may be a subject of future work. Another question which goes beyond the scope of the present paper is whether continuous hysteresis relationships could be thought of as some continuous limit of the gradient systems.

The paper is organized as follows. In the next section we discuss the setting of the problem and 
present the main result. In Section 3, a graph describing the dynamics of the discrete Preisach nonlinearity is considered as an example. Section 4 contains the proofs. We use an explicit construction for the functional of the gradient system in order to obtain a realization of any given graph that encodes transitions between the states of a hysteretic system. Finally, a brief discussion of multi-scale dynamics and their adiabatic
limit is presented in the Apprndix.

\section{Main result}
The problem can be formally stated as follows. Consider a directed graph $\Gamma$
with the following properties.

\medskip
(i) The set $S$ of all vertices is a union of $n+1$ non-empty disjoint subsets:
$$
S=\cup_{i=0}^n S_i;\qquad S_i\cap S_j=\emptyset, \ i\ne j;\qquad S_i\ne\emptyset.
$$

\medskip
(ii) For every $1\le i\le n-1$ and every vertex $s\in S_i$ there are exactly two edges $s s'$ and $s s''$ emanating from $s$.
The end point $s'$ of one edge belongs to $S_{i-1}$, the end point $s''$ of the other edge belongs to $S_{i+1}$.
Furthermore, for every $s\in S_0\cup S_{n}$ there is exactly one edge $ss'$ emanating from $s$; the end point $s'$
of this edge belongs to $S_1$ if $s\in S_0$ and to $S_{n-1}$ if $s\in S_n$.

\medskip
Graphs with these two properties will be called {\em admissible}.

Dynamics on an admissible  graph $\Gamma$ is defined as follows.
It is assumed that the input $u$ of the hysteretic system takes values from a finite set
$$
U=\{u^0,\ldots,u^n\},\qquad u^0< u^1<\ldots <u^n.
$$
Any input sequence $u_t$ should satisfy $u_t\in U$ and either $u_{t}=u^{i-1}$ or $u_t=u^i$ or $u_t=u^{i-1}$ if $u_t=u^i$ for all $t$.
Vertices $s$ from the set $S_i$ represent possible states of the hysteretic system at any moment when $u_t=u^i$.
Edges of the graph define transitions between the states induced by the input sequence.
If $u_{t-1}=u^i$, $u_t=u^{i+ \sigma}$ with $\sigma=\pm1$, and $s_{t-1}=s'\in S_i$, then $s_t=s''$ where $s's''$ is a unique directed edge with $s''\in S_{i+\sigma}$.

Now, let us consider the gradient system
\begin{equation}\label{3}
\dot {\bm x} = -\nabla_x V(\bm{x},u)
\end{equation}
with $\bm{x}\in \mathbb{R}^d,\ u\in \mathbb{R}$.
Question is whether for any given admissible graph $\Gamma$ it is possible to construct a functional $V(\bm{x},u)$
such that the anti-gradient dynamics defined by $V$ agrees with the dynamics on $\Gamma$. A few more definitions are in order to formalize this question.
We will assume everywhere that $V$ is continuously differentiable in $\bm{x}$ and continuous with respect to $\bm{x},u$.

Suppose that the functional $V(\cdot,u)$ has an isolated critical point
$\bm{x}^*(u)$ for each $u\in [u_-,u_+]$ and that this point continuously depends on $u$.
Suppose that $\bm{x}^{*}(u)$ is a point of local minimum for $u\in [u_-,u_+)$, which
undergoes a saddle-node bifurcation at
$u=u_+$. That is, $\bm{x}^*(u_+)$ is a saddle-node point of the gradient field $-\nabla_x V(\cdot;u_*)$
and a ball $|\bm{x}-\bm{x}^*(u_+)|<\rho$ does not contain critical points of $V(\cdot,u)$ for some interval $u\in [u_+,u_++\delta)$.
Assume that $\bm{x}^{**}(u_+)$ is an isolated local minimum point of $V(\cdot,u_+)$.
We say that the local minimum point $\bm{x}^*(u_-)$ of $V(\cdot,u_-)$
{\em connects to} the local minimum  point $\bm{x}^{**}(u_+)$ of $V(\cdot,u_+)$ if
the unstable manifold of the saddle-node equilibrium point $\bm{x}^*(u_+)$ of the gradient field $-\nabla_x V(\cdot,u_+)$ belongs to the basin of attraction
of the stable equilibrium $\bm{x}^{**}(u_+)$ of this field.

If $u$ changes adiabatically (very slowly) from $u_-$ to $u_+$
and gradient system \eqref{3} is at the equilibrium $\bm{x}^*(u_-)$ for $u=u_-$,
then this system will remain at the equilibrium $\bm{x}^*(u)$ until the adiabatic input $u$
crosses the value $u=u_+$. At this point, the equilibrium $\bm{x}^*(u)$ will disappear in the saddle-node bifurcation
and the system will switch to another equilibrium $\bm{x}^{**}(u_+)$ along the unstable manifold
of $\bm{x}^*(u_+)$. Some brief discussion of multi-scale systems that lead to equation \eqref{3}
in the adiabatic limit is presented in the Appendix.

Similarly to the above definition, we say that a local minimum point $\bm{x}^*(u_+)$ of $V(\cdot,u_+)$
{\em connects to} a local minimum  point $\bm{x}^{**}(u_-)$ of $V(\cdot,u_-)$ if
$\bm{x}^*(u)$ undergoes the saddle node-bifurcation at $u=u_-$ and
the unstable manifold of $\bm{x}^*(u_-)$ belongs to the basin of attraction
of the stable equilibrium $\bm{x}^{**}(u_-)$ of the gradient flow.

We will also say that $\bm{x}^*(u_-)$ and $\bm{x}^{*}(u_+)$ are {\em reversibly connected}
if $\bm{x}^*(u)$ is an isolated local minimum point of $V(\cdot,u)$ for all $u\in[u_-,u_+]$
and $\bm{x}^*(u)$ continuously depends on $u$. For such equilibrium points,
the gradient system remains at the equilibrium state $\bm{x}^*(u)$ at all times as long as the input is varied adiabatically
within the range $[u_-,u_+]$ and provided that the system was at this state initially.

Let us extend the above definitions to the case of several transitions between the equilibrium points points.
Suppose that a local minimum point $\bm{x}^*_{k-1}$ of the functional $V(\cdot,v_{k-1})$ connects
to the local minimum point $\bm{x}_{k}^*$ of the functional $V(\cdot,v_k)$ for every $k=1,\ldots,\ell$, where
$u_-=v_0<v_1<\cdots<v_\ell=u_+$ is a partition  of the interval $[u_-,u_+]$.
Then we say that $\bm{x}_0^*$ {\em connects to} $\bm{x}_\ell^*$. Similarly,
we say that $\bm{x}_\ell^*$ {\em connects to} $\bm{x}_0^*$ if $\bm{x}_{k}^*$ connects to $\bm{x}^*_{k-1}$ for every $k=1,\ldots,\ell$.

We call a gradient system $\dot {\bm{x}}=-\nabla_x V(\bm{x},u)$, $\bm{x}\in \mathbb{R}^d$ a {\em realization} of an admissible graph $\Gamma$ if

\medskip
(i) The functional $V(\cdot,u)$ is radially increasing for each $u\in [u_0,u_n]$:
$$
\lim_{|\bm{x}|\to \infty}V(\bm{x},u)=\infty.
$$

\medskip
(ii) The set of the points of local minimum of the functional $V(\cdot,u_i)$
is in one-to-one correspondence $\bm{x}^*=X_i(s)$ with the subset $S_i\ni s$ of the vertices of $\Gamma$ for every $i=0,\ldots,n$.

\medskip
(iii) A local minimum point $\bm{x}^*$ of the functional $V(\cdot,u_i)$ connects to
a local minimum point $\bm{x}^{**}$ of the functional $V(\cdot,u_{i+\sigma})$ with $\sigma=\pm1$
iff there is an edge $s's''$ of $\Gamma$ such that $\bm{x}^*=X_i(s')$, $\bm{x}^{**}=X_{i+\sigma}(s'')$.

\medskip
Adiabatic variation of the input $u$ of the gradient system over the interval $[u^0,u^n]$
results in transitions from one to another equilibrium state.
If we restrict observations to the instants when the input value belongs to the finite set $U=\{u^0,\ldots, u^n\}$, then
according to the above definition,
the dynamics of transitions between the equilibrium states of the gradient system realizing the graph $\Gamma$ will be identical to dynamics on the graph $\Gamma$.

The following statement is the main result of this paper.

\begin{theorem}\label{theorem1}
There is a planar realization for every admissible graph $\Gamma$.
\end{theorem}

A constructive proof of this theorem is presented in Section \ref{proof}.

\section{Example: Realization of the Preisach hysteresis model}
Fig.~\ref{Fig1} defines one of the simplest hysteresis operators called the non-ideal relay.
The state of the relay corresponding to the right minimum of the double well potential is
denoted by $0$, the other state corresponding to the left minimum is denoted by $1$.
The relay switches from state $0$ to state $1$ when the input $u=u_t$ exceeds a threshold value
$u_+$ and from state $1$ to state $0$ when the input becomes smaller than
a different threshold value $u_-$ which satisfies $u_-<u_+$.
Therefore $[u_-,u_+]$ is a bi-stability interval for the relay, see Fig.~\ref{Fig2} (left).
For a given input sequence $u_t$ we will denote the time series of the varying state
of the non-ideal relay 
by
$$
r_t={\mathcal R}_{u_{-},u_+}[u_t].
$$

\begin{figure}
    \centering
        \includegraphics[width=0.4\textwidth]{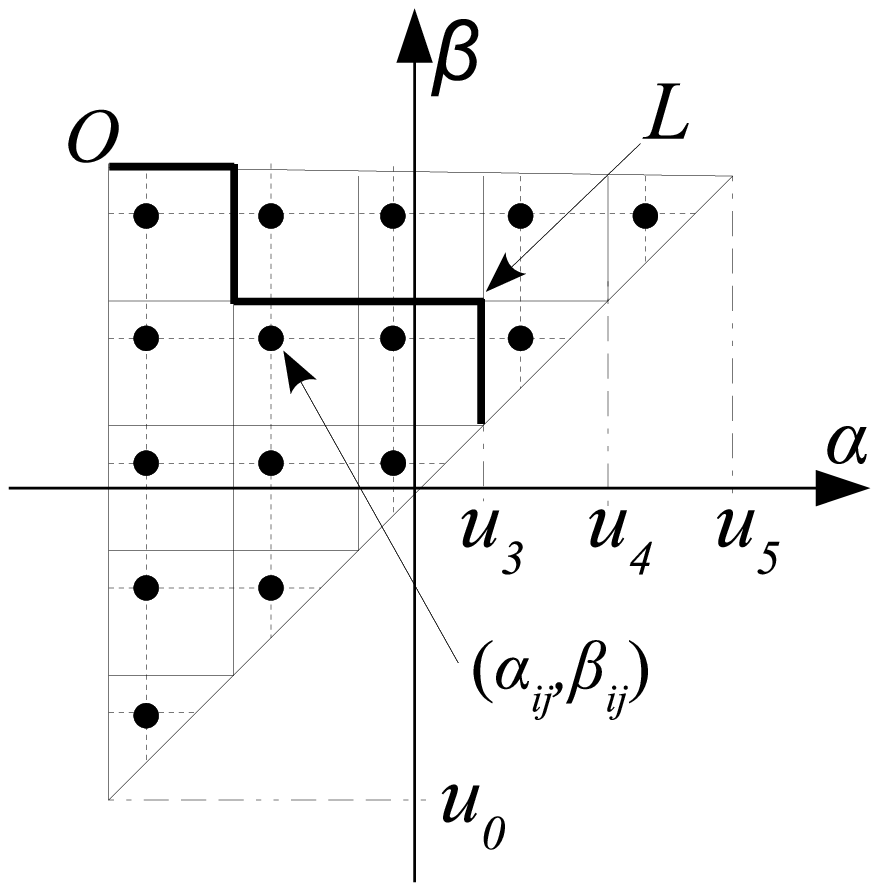} \qquad \includegraphics[width=0.45\textwidth]{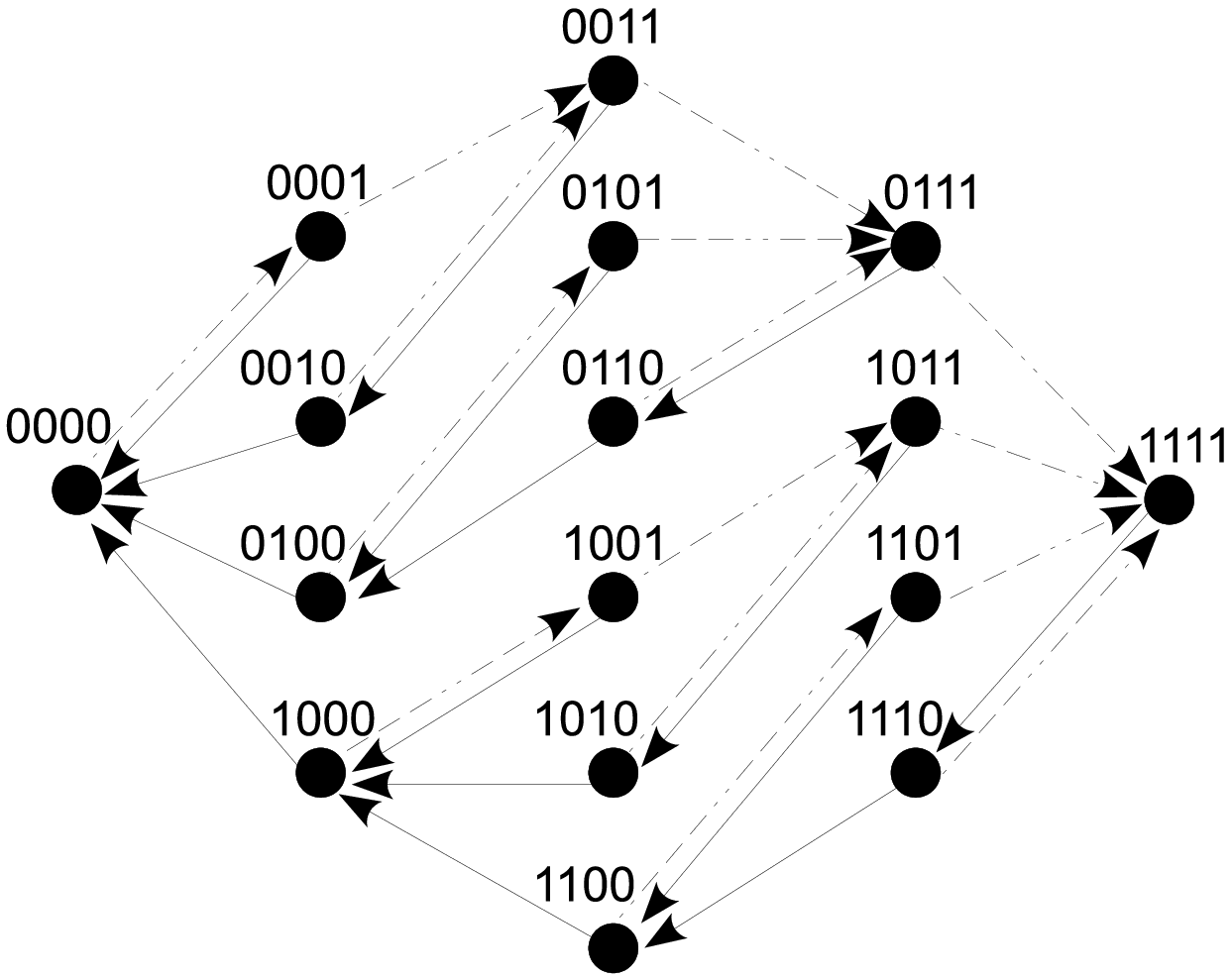}
    \caption{{\bf Left:} $(\alpha,\beta)$ plane for the discrete Preisach model with $N(N+1)/2=15$ relays ($N=5$). Each black point represents
    a non-ideal relay with thresholds $\alpha,\beta$. The staircase line $L$ represents a state of the model. Those relay below this line are in state $1$,
    the relays above this line are in state $0$.
    {\bf Right:} Graph $\Gamma$ describing the discrete Preisach model for $N=4$. The graph has $2^N=16$ nodes corresponding to
    $16$ possible staircase states
    of the model, each encoded by a four-tuple of zeros and ones. The sets $S_0,S_1,S_2,S_3,S_4$ consist of $1,4,6,4,1$ nodes, respectively.} \label{Fig3}
\end{figure}

The finite (discretized) Preisach model is a superposition of $N(N+1)/2$ relays that have a common input $u=u_t$ which often is also discretized. The relays have different thresholds.
The output of the model is defined by the formula
$$
p_t=\sum_{i=1}^{N}\sum_{j=1}^{i} {\mathcal R}_{\alpha_{ij},\beta_{ij}}[u_t]
$$
where the thresholds of the relays, $(\alpha_{ij},\beta_{ij})$, typically sit on a lattice of step $\varepsilon$ in the half-plane $\beta\ge \alpha$ of a plane $(\alpha,\beta)$
and form a mesh in an isosceles right triangle;
and, the input takes values $\{u^0,\ldots,u^N\}$ at the nodes situated between the threshold values, see Fig.~\ref{Fig3} (left).
States of the model can be identified with the polylines $L$ that connect the vertex $O$ at right angle of the triangle with its hypothenuse on the line $\alpha=\beta$, see Fig.~\ref{Fig3} (left).
Each polyline has $N$ links that have length $\varepsilon$ and are either horizontal or vertical. Therefore a state can be encoded by, and identified with, an $N$-tuple $(a_1,\ldots,a_N)$
where each $a_i$ equals either 0 or 1 and we write 0 for a vertical link and 1 for a horizontal link; $a_1$ defines the direction (vertical or horizontal) for the link starting at the point $O$,
$a_2$ correspond to the next link and so on. The set $S_i$ corresponding to the input value $u^i$ consists of the states $L$ that have exactly $i$ entries 1.

A detailed description of the Preisach model can be found, for example, in \cite{BroSpr, Pre, Pre1, Pre2, Pr, Pre3}.
In particular, transitions between the states are defined by simple rules.
Adapting these rules for the discretized model considered here,
we obtain the following description of transitions between the states.
If the system is at a state $(a_1,\ldots,a_N)\in S_i$ and the input changes from $u^i$ to $u^{i+1}$,
then the system transits to the state $(b_1,\ldots,b_N)$ defined by the formulas
$$
b_k=1 \quad {\rm if} \quad a_{k}=0 \ {\rm and\ either} \ k=N \ {\rm or} \ a_{k+1}=\cdots=a_N=1;\quad
b_j=a_j \ {\rm for} \ j\ne k.
$$
That is, the last 0 in the $N$-tuple $(a_1,\ldots,a_N)$ changes to 1.
Similarly, when the input changes from $u^i$ to $u^{i-1}$,
the system transits from a state $(a_1,\ldots,a_N)\in S_i$ to the state $(b_1,\ldots,b_N)$ defined by
$$
b_k=0 \quad {\rm if} \quad a_{k}=1 \ {\rm and\ either} \ k=N \ {\rm or} \ a_{k+1}=\cdots=a_N=0;\quad
b_j=a_j \ {\rm for} \ j\ne k,
$$
that is, the last 1 in the $N$-tuple $(a_1,\ldots,a_N)$ changes to 0.
These rules define the edges of the graph $\Gamma$. Fig.~\ref{Fig3} (right)
shows $\Gamma$ for $N=4$.

\section{Proofs}\label{proof}

\subsection{Proof of Theorem \ref{theorem1}}

\begin{figure}
    \centering
        \includegraphics[width=0.4\textwidth]{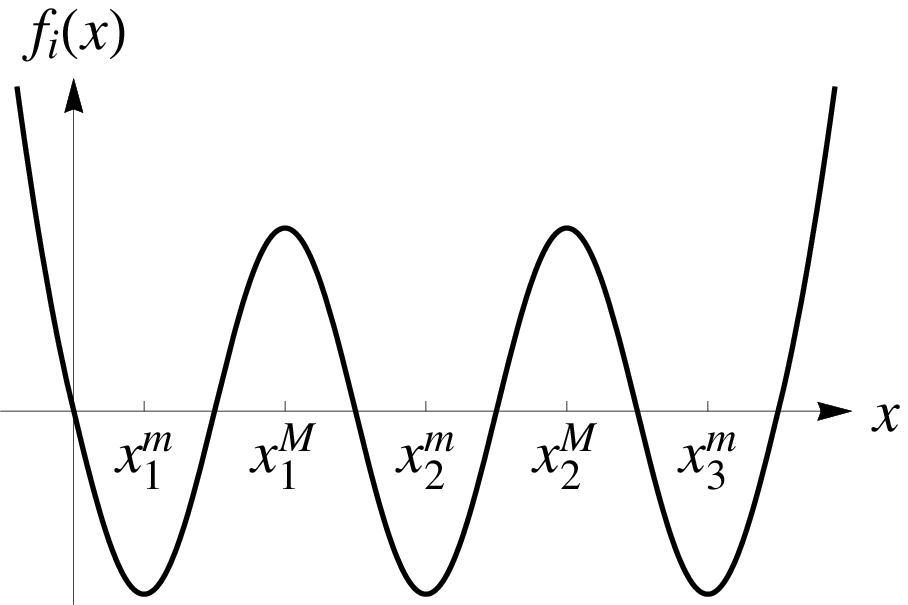} \qquad\quad \includegraphics[width=0.45\textwidth]{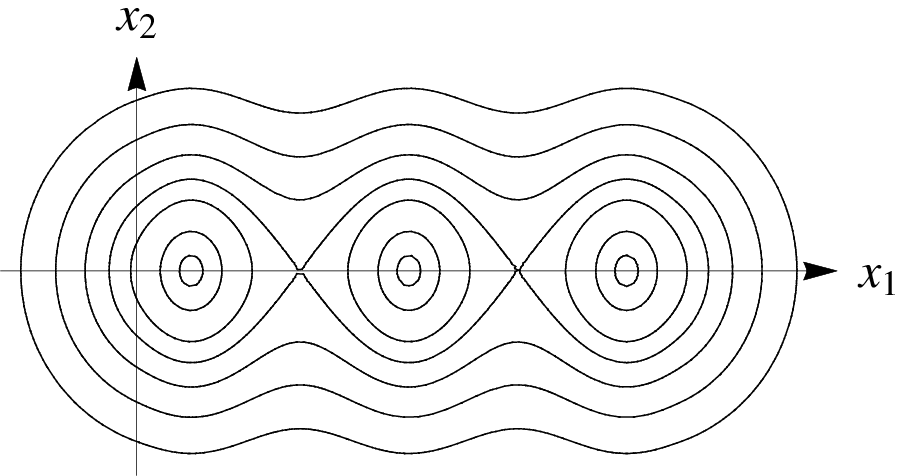}
    \caption{ Function $f_i:\mathbb{R}\to \mathbb{R}$ (left) and level lines of the functional $V(\bm{x},u^i)=f_i(x_1)+x_2^2$ (right).} \label{Fig4}
\end{figure}

Realization of the system for $u=u^i$ will be a smooth functional
$$
V(\bm{x},u^i)=f_i(x_1)+x_2^2, \qquad \bm{x}=(x_1,x_2)\in \mathbb{R}^2,
$$
which has $N_i$ local minimum points $\bm{x}_1^m=(x^m_1,0),\ldots,\bm{x}_{N_i}^m=(x^m_{N_i},0)$, where $x^m_1<\cdots<x^m_{N_i}$ are the local minima of the scalar function $f_i: \mathbb{R}\to \mathbb{R}$; $N_i$ is a number of vertices in the set $S_i$;
and, $f_i(x_1)=y_1^2$ for sufficiently large $|x_1|$, see Fig.~\ref{Fig4}.
We require for convenience that $f_i(x^m_1)=\cdots=f_i(x^m_{N_i})=m_*$ as well as
$f_i(x^M_1)=\cdots=f_i(x^M_{N_i-1})=M_*>m_*$ for the local maximum points of $f_i$ that separate its local minima (although this is not necessary for the following construction),
see Fig.~\ref{Fig4}.

For the intermediate value $u^{i+1/2}=(u_i+u_{i+1})/2$, let us define the functional $V(\cdot,u^{i+1/2})$ in a similar fashion so that
it would have $N_i+N_{i+1}$ minimum points. That is, $V(\cdot,u^{i+1/2})=f_{i+1/2}(x_1)+x_2^2$ where
$f_{i+1/2}(x^m_1)=\cdots=f_{i+1/2}(x^m_{N_i})=f_{i+1/2}(\tilde x^m_1)=\cdots=f_{i+1/2}(\tilde x^m_{N_{i+1}})=m_*$ at the local minimum points
and $f_{i+1/2}(x^M_1)=\cdots=f_{i+1/2}(x^M_{N_i})=f_{i+1/2}(\tilde x^M_1)=\cdots=f_{i+1/2}(\tilde x^M_{N_{i+1}-1})=M_*$ at the local maximum points of $f_{i+1/2}$ as well as
$f_{i+1/2}(x_1)=x_1^2$ for sufficiently large $|x_1|$.

The main step of the proof is the following lemma.

\begin{lemma}\label{lemma1}
Suppose that a smooth functional $V_0(\bm{x})=V_0(x_1,x_2)=f(x_1)+x_2^2$ has $N$ isolated local minimum points
and $f(x_1)=x_1^2$ for sufficiently large $|x_1|$. Let $(j_1,\ldots,j_N)$ be a permutation of numbers $(1,\ldots,N)$.
Then there is a smooth deformation $V(\cdot,u)$, $u\in[u_-,u_+]$, such that
$$
V(\cdot,u_-)=V(\cdot,u_+)=V_0(\cdot)
$$
and the functional $V(\cdot,u)$ has $N$ isolated local minimum points $\bm{x}_i^*(u)$
for each  $u\in[u_-,u_+]$. Moreover, the continuous functions
$\bm{x}_i^*(u): [u_-,u_+]\to \mathbb{R}^2$ satisfy
$$
\bm{x}_k^*(u_-)=\bm{x}_{j_k}^*(u_+),\quad \quad
k=1,\ldots,N.
$$
\end{lemma}

A constructive proof of this statement is presented in the next section.

We will say that the deformation defined in Lemma \ref{lemma1} {\em permutes} the minimum points.
Indeed, for every $k$, the minimum point $\bm{x}_k^*(u_-)$ of $V(\cdot,u_-)$ and the minimum point $\bm{x}_{j_k}^*(u_+)$ of $V(\cdot,u_+)$ are reversibly connected.

Assume that a smooth function $f:\mathbb{R}\to \mathbb{R}$ that has $N$ local minimum points $x^m_k$ separated by $N-1$ local maximum points $x^M_k$,
\begin{equation}\label{n1}
x_1^m<x_1^M<x_2^m<x_2^M<\dots<x_{N-1}^m<x_{N-1}^M<x_N^m,
\end{equation}
and grows at infinity as in Lemma \ref{lemma1}. Assume that
a smooth function $\tilde f: \mathbb{R}\to \mathbb{R}$ has $N-1$ local minimum points $\tilde x^m_k$ separated by points of local maximum $\tilde x^M_k$,
\begin{equation}\label{n2}
\tilde x_1^m<\tilde x_1^M<\tilde x_2^m<\tilde x_2^M<\dots<\tilde x_{N-2}^M<\tilde x_{N-1}^m,
\end{equation}
and grows at infinity in the same way. Consider a smooth deformation $f(\cdot,u)$, $u\in[u_-,u_+]$, with $f(\cdot,u_-)=f(\cdot)$, $ f(\cdot,u_+)=\tilde f(\cdot)$
such that the maximum point  $x_{N-1}^M$  and the minimum point $x_N^m$ collide and disappear in the saddle node bifurcation
for some $u_{SN}\in (u_-,u_+)$, while the other extremum points continue from $u_-$ to $u_+$, see Fig.~\ref{Fig1}. That is,
$f(\cdot,u)$ has $2N-1$ critical points for $u\in[u_-,u_{SN})$ and $2N-3$ critical points for $u\in(u_{SN},u_+]$, each
critical point $x_k^m=x_k^m(u)$, $x_k^M=x_k^M(u)$ is isolated, and they satisfy relations \eqref{n1} for $u\in[u_-,u_{SN})$,
relations \eqref{n2} for $u\in(u_{SN},u_+]$ and the relation $x_{N-1}^M(u_{SN})=x_N^m(u_{SN})$ at the bifurcation point $u=u_{SN}$.
Furthermore, all the critical points $x_1^m(u),x_1^M(u),\ldots,x_{N-2}^M(u), x_{N-1}^m(u)$ continuously depend on $u$ on the interval $[u_-,u_+]$
and the critical points $x_{N-1}^M(u), x_N^m(u)$ continuously depend on $u$ in their domain $[u_-,u_{SN}]$.
We will say that a deformation satisfying these conditions {\em connects} the local minimum $x_N^m=x_N^m(u_-)$ of $f$
to the local minimum $\tilde x_{N-1}^m=x_{N-1}^m(u_+)$ of $\tilde f$. Without loss of generality, we
will impose an extra condition $\partial V/\partial u (\cdot,u_-)=\partial V/\partial u (\cdot,u_+)=0$.
Such deformations can be easily constructed explicitly.

The above deformation $f(\cdot,u)$ induces the homotopy
$$
V(\bm{x},u)=f(x_1,u)+x_2^2,\qquad \bm{x}=(x_1,x_2)\in \mathbb{R}^2,\ u\in[u_-,u_*],
$$
that connects the local minimum point $\bm{x}_N^*(u_-)=(x_N^m(u_-),0)$ of the functional $V(\cdot,u_-)$
to the local minimum  $\bm{x}_{N-1}^*(u_-)=(x_{N-1}^m(u_+),0)$ of the functional $V(\cdot,u_+)$.
Every other local minimum $\bm{x}_k^*(u_-)=(x_k^m(u_-),0)$ of $V(\cdot,u_-)$
is reversibly connected with the local minimum $\bm{x}_k^*(u_+)=(x_k^m(u_+),0)$ of $V(\cdot,u_+)$ ($k=1,\ldots,N-2$)
along this deformation.

A deformation $V(\cdot,u)$, $u\in[u_-,u_+]$,  will be called a {\em concatenation} of deformations
$V_k(\cdot,u)$, $u\in[v^k,v^{k+1}]$ if $V(\cdot,u)=V_k(\cdot,u)$ for $u\in[v^k,v^{k+1}]$,
where $u_-=v^0<v^1<\cdots<v^{K-1}<v^K=u_+$ is a partition of the interval $[u_-,u_+]$.

We now define the function $V(\cdot,\cdot)$ for
$\bm{x}\in\mathbb{R}^2$, $u\in [u^i,u^{i+1}]$ as a concatenation of deformations $V(\cdot,u)$ that permute minima
and deformations that eliminate minima via a saddle-node bifurcation.
Assume that the set $S_i$ of the graph $\Gamma$ has $N$ vertices $s_k$
and the set $S_{i+1}$ has $N'$ vertices $s_k'$.
Using the construction described above, we define the functionals
$V(\cdot, u^i)$, $V(\cdot, u^{i+1})$, $V(\cdot, u^{i+1/2})$ with $N$, $N'$ and $N+N'$
minimum points, respectively. Now, consider an arbitrary partition of the segment
$[u^{i+1/2},u^{i+1}]$ into $2N+1$ intervals, $u^{i+1/2}=v_0<v_1<\cdots<v_{2N+1}=u^{i+1}$,
and a partition of the segment $[u^{i},u^{i+1/2}]$ into $2N'+1$ intervals,
$u^{i}=w_{2N'+1}<\cdots<w_1<w_0=u^{i+1/2}$.
Let us associate a local minimum $\bm{x}_{k}^*$ of the functional $V(\cdot,u^i)$ with a vertex
$s_k\in S_i$ of $\Gamma$ (as one-one correspondence),
each minimum $\bm{x}_{j}^{**}$ of $V(\cdot,u^{i+1})$ with a vertex $s_j'\in S_{i+1}$,
and each minimum $\bm{x}^\dagger_{\ell}$ of $V(\cdot,u^{i+1/2})$ with one point of $S_i\cup S_{i+1}$.
Next, let us associate the point $v_{2k}$, $k=1,\ldots,N$, with the
edge $s_k s_{j_k}'$ of $\Gamma$ that originates at the vertex $s_k\in S_i$ and the points $w_{2k}$ with the edge $s_k's_{\ell_k}$ where $s_k'\in S_{i+1}$.
Consider a deformation $V(\cdot,u)$, $u\in [v_0,v_1]$, that permutes
the local minima and connects the minimum $\bm{x}^\dagger_1$ of $V(\cdot,v_0)$
associated with the vertex $s_1$ with the largest local minimum of $V(\cdot,v_1)$
(with respect to the ordering defined by the $x_1$-coordinate of the minimum)
and the minimum $\bm{x}^\dagger_{N+j_1}$ of $V(\cdot,v_0)$
associated with the vertex $s_{j_1}'$ with the second largest local minimum of $V(\cdot,v_1)$.
Then consider a deformation $V(\cdot,u)$, $u\in [v_1,v_2]$, that eliminates the largest
minimum of $V(\cdot,v_1)$ via the saddle-node bifurcation.
A concatenation of these two deformations reversibly connects the minimum point $\bm{x}^\dagger_{N+j_1}$ of  $V(\cdot,v_0)$
associated with $s_{j_1}'$ with the largest minimum point of  $V(\cdot,v_2)$
as well as connects the minimum point $\bm{x}^\dagger_{1}$ of  $V(\cdot,v_0)$
associated with $s_1$ to the same largest minimum point of  $V(\cdot,v_2)$.
Thus, this concatenation realizes the transition $s_1s_{j_1}'$ as $u$ changes from $v_0=u^{i+1/2}$
to $v_2$. Furthermore, every other minimum of  $V(\cdot,v_0)$ is reversibly connected with
one and only one minimum of  $V(\cdot,v_2)$. The total number of minima
of $V(\cdot,v_2)$ is $N+N'-1$, one less than $V(\cdot,v_0)$ has.

Now, we proceed by induction to extend the deformation from an interval $[v_0,v_{2k-2}]$ to the
interval $[v_0,v_{2k}]$. The functional $V(\cdot,v_{2k-2})$ has a minimum $\bm{x}^\ddagger_k$
reversibly connected with the minimum $\bm{x}_k^\dagger$ of $V(\cdot,v_0)$, which is associated with the vertex $s_k$. Also, $V(\cdot,v_{2k-2})$
has a minimum $\bm{x}^\ddagger_{N+j_k}$ reversibly connected with the minimum $\bm{x}_{N+j_k}^\dagger$ of $V(\cdot,v_0)$, which is associated with the vertex $s_{j_k}'$ (recall that $v_{2k}$ is associated with the edge $s_k s_{j_k}'$ of $\Gamma$).
We use a deformation $V(\cdot,u)$, $u\in[v_{2k-2},v_{2k-1}]$, that permutes the minima
and connects the minimum points $\bm{x}^\ddagger_k$ and  $\bm{x}^\ddagger_{N+j_k}$ of $V(\cdot,v_{2k-2})$ with the largest
minimum and the second largest minimum of the functional $V(\cdot,v_{2k-1})$, respectively.
We further use a deformation $V(\cdot,u)$, $u\in[v_{2k-1},v_{2k}]$ that eliminates the largest
minimum. A concatenation of these deformations with the deformation defined
on the interval $v_0\le u\le v_{2k-2}$ completes the definition of $V$ on the interval
$v_0\le u\le v_{2k-2}$ and the inductive step. After $N$ steps, we obtain the deformation
defined for $u\in [v_0,v_{2N}]$. Finally, on the last interval $v_{2N}\le u\le v_{2N+1}=u^{i+1}$
we can use any deformation $V(\cdot,u)$ that permutes $N'$ local minima $\bm{z}_1^*,\ldots, \bm{z}_{N'}^*$
of $V(\cdot,v_{2N})$ and reversibly connects a local minimum $\bm{x}^{**}_{k}$ of $V(\cdot,v_{2N+1})$
with those (unique) local minimum $\bm{z}^*_j$ of $V(\cdot,v_{2N})$ that is reversibly connected with the minimum
$\bm{x}^\dagger_{N+k}$ of $V(\cdot,v_{0})$ (for every $k=1,\ldots,N'$). Now, the concatenation of this deformation with the deformation $V(\cdot,u)$, $u\in[v_0,v_{2N}]$,
realizes all the edges $s_ks_{j_k}'$ of $\Gamma$ for the input increasing from $u^{i+1/2}$ to $u^{i+1}$ and, simultaneously,
establishes the reversible connection between the minimum points $\bm{x}^\dagger_{N+k}$ and $\bm{x}^{**}_k$ of $V(\cdot,u^{i+1/2})$ and $V(\cdot,u^{i+1})$, respectively, for each $k=1,\ldots,N'$.

Using the same construction, we define a deformation $V(\cdot,\cdot)$ on the interval $[u^i,u^{i+1/2}]$
that realizes all the edges $s_k' s_{\ell_k}$ of $\Gamma$ with $s_k'\in S_{i+1}$, $s_{\ell_k}\in S_i$ for the input decreasing from $u^{i+1/2}$ to $u^{i}$ and, simultaneously,
establishes the reversible connection between the minimum points $\bm{x}^\dagger_{j}$ and $\bm{x}^{*}_j$ of $V(\cdot,u^{i+1/2})$ and $V(\cdot,u^{i})$, respectively, for each $j=1,\ldots,N$.
By construction, the concatenation of the two deformations $V(\cdot,\cdot)$ defined on the intervals $[u^i,u^{i+1/2}]$ and $[u^{i+1/2},u^{i+1}]$
realizes all the edges $s_k s_{j_k}'$ with $s_k\in S_i$, $s_{j_k}'\in S_{i+1}$ and all the edges $s_k's_{\ell_k}$ with $s_k'\in S_{i+1}$, $s_{\ell_k}\in S_i$ of the graph $\Gamma$.
Finally, the concatenation of such deformations defined on all the intervals $[u^i,u^{i+1}]$ realizes the whole graph $\Gamma$.
This completes the proof of Theorem \ref{theorem1}.

\subsection{Proof of Lemma \ref{lemma1}}
Suppose that the functional $V_0(\bm{x})=f(x_1)+x_2^2$ has exactly $N$ local minimum points
$\bm{x}_1^*=({x}_1^m,0),\ldots, \bm{x}=({x}_N^m,0)$. Here the isolated local minimum points $x_i^m$ of the scalar function $f$ are
separated by the isolated local maximum points
 $$x_1^m<x_1^M<x_2^m<\cdots<x_{N-1}^M<x_N^m$$
 and we assume that $f(x_1^m)=\cdots=f(x_N^m)=m_*<M_*=f(x_1^M)=\cdots=f(x_{N-1}^M)$.
It suffices to prove the lemma for a permutation that interchanges the points $\bm{x}_j^*$ and $\bm{x}_{j+1}^*$ because
any permutation can be obtained as a finite sequence of such ``elementary''  permutations and therefore a realization
of any permutation can be obtained as a concatenation of the realizations of the ``elementary'' permutations
each interchanging two neighboring minimum points.

Starting with the functional
$
V(\bm{x},u_-)=V_0(\bm{x})
$
we will concatenate several deformations to obtain the functional $V(\bm{x},u_+)=V(\bm{x},u_-)$ with $\bm{x}^*_j(u_+)= \bm{x}^*_{j+1}(u_-)$, $\bm{x}^*_{j+1}(u_+)= \bm{x}^*_{j}(u_-)$ and $\bm{x}^*_i(u_+)= \bm{x}^*_{i}(u_-)$ for $i\ne j, j+1$.

\begin{figure}
    \centering
        \includegraphics[width=0.8\textwidth]{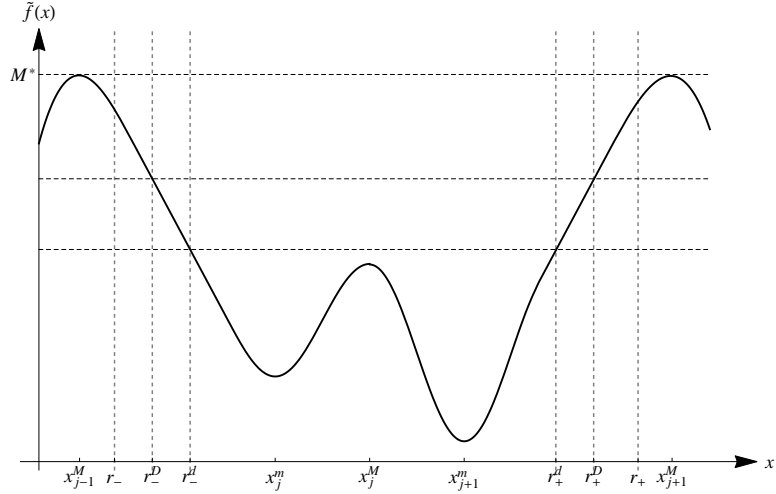} 
    \caption{Function $\tilde f$ and points  $x_{j-1}^M<r_-<r_-^D<r^d_-<x_j^m<x_j^M<x_{j+1}^m<r_+^d<r_+^D<r_+<x^M_{j+1}$.} \label{Fig5}
\end{figure}

First, define a smooth function $\tilde f(x_1)$ such that $\tilde f$ and $f$ have the same extremum points $x_i^m, x_i^M$;
the maximum value of $\tilde f$ at the point $x_j^M$ satisfies $\tilde f(x_j^M)<M_*$ where $M_*$ is the value of $f$ at all its maximum points;
there are points $r_-,r_+$ satisfying $x_{j-1}^M<r_-<x_j^m<x_j^M<x_{j+1}^m<r_+<x^M_{j+1}$ such that $\tilde f(x_1)=f(x_1)$ for $x_1\le r_- $ and $x_1\ge r_+$
(we agree that $x_{0}^M$, $x_{N}^M$ are the non-critical points at which $f$ equals $M_*$)
and $\tilde f(r_-)>\tilde f (x_j^M)$, $\tilde f(r_+)>\tilde f (x_j^M)$;
and, there are two intervals $[r_-^D,r_-^d]\subset (r_-,x_j^m)$
and $[r_+^d,r_+^D]\subset (x_{j+1}^m,r_+)$ such that $r^d_- -r^D_-=r^D_+-r^d_+$ and
\begin{equation}\label{sym}
\tilde f(r_-^D+s)=\tilde f(r_+^D-s),\quad \ \  s\in[0,r^d_- -r^D_-],
\end{equation}
\begin{equation}\label{condi}
\tilde f(r_\pm^d)>\tilde f (x_j^M),
\end{equation}
see Fig.~\ref{Fig5}.
Without loss of generality we can assume that
$$
|x_j^m-x_j^M|=|x_{j+1}^m-x_j^M|,\quad |r_-^d-x_j^M|=|r_+^d-x_j^M|,\quad |r_-^D-x_j^M|=|r_+^D-x_j^M|.
$$
The first deformation we use is the linear deformation
$$
V(x_1,x_2,u)=(1-a(u)) f(x_1)+a(u) \tilde f(x_1)+x_2^2,\quad u\in[u_-,v^1],
$$
with $v^1<u_+$, where $a(u_-)=0$, $a(v^1)=1$ and $\frac{d a}{du}(u_-)=\frac{d a}{du}(v^1)=0$
(the latter condition will be used to smoothly concatenate this deformation with the next one).

\begin{figure}
    \centering
     \begin{minipage}{0.44\textwidth}   \includegraphics[width=1\textwidth]{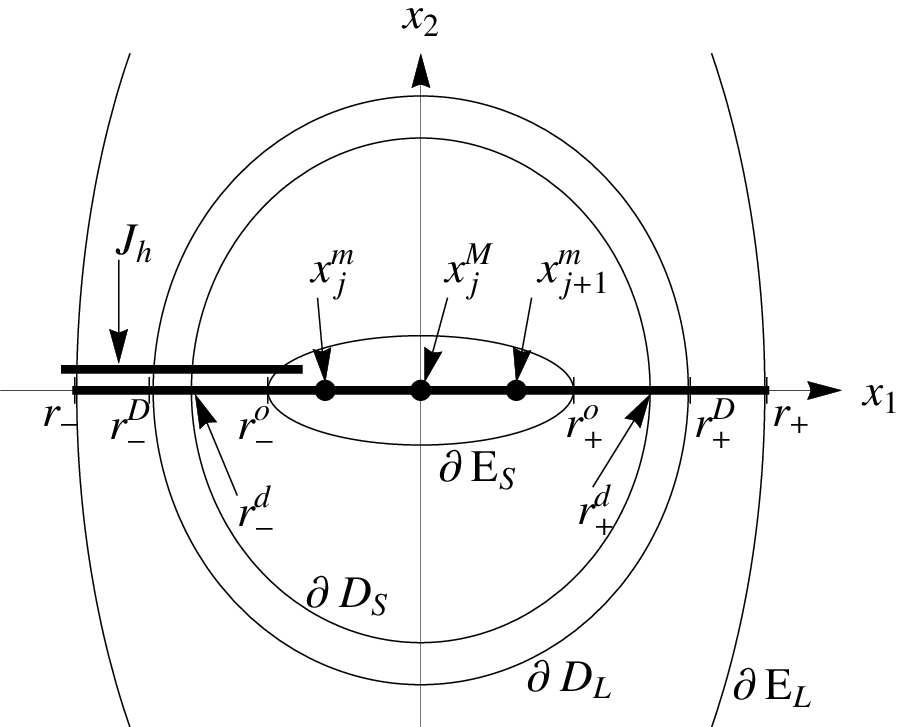}\end{minipage} \qquad \begin{minipage}{0.43\textwidth} \includegraphics[width=1\textwidth]{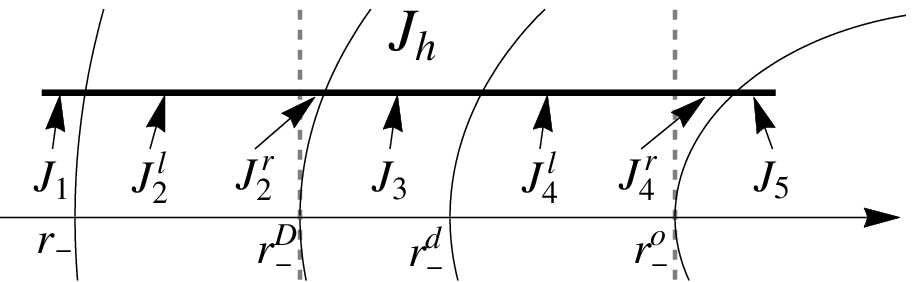}\end{minipage}
    \caption{{\bf Left:} The functionals $V(\cdot, v^1)$ and $\Phi(\cdot)$ coincide inside the ellipse $E_S$, outside the ellipse $E_L$ and on the segment $J$. The level lines of the functional $\Phi(\cdot)$
    inside the annulus $D_L\setminus D_S$ are concentric circles with the center at the point $(x^M_j,0)$. Therefore deformation \eqref{rotate} is continuous. {\bf Right:} Zoom of the segment $J_h$ and its subdivision
    into small segments by the lines $\partial E_L, \partial D_L, \partial D_S, \partial E_S$ and $x_2=r^D_-$, $x_2=r^0_-$.} \label{Fig6}
\end{figure}

Next, we introduce two more points $r^0_\pm$ such that
$$
x_{j-1}^M<r_-<r_-^D<r_-^d<r_-^0< x_j^m<x_j^M<x_{j+1}^m< r_+^0<r_+^d<r_+^D<r_+<x^M_{j+1}.
$$
Let us define two concentric circles $\partial D_L$ and $\partial D_S$ with the center $(x_j^M,0)$
and the radii $(r_+^D-r_-^D)/2=x_j^M-r_-^D=r^D_+-x_j^M$ and $(r_+^d-r_-^d)/2=x_j^M-r_-^d=r^d_+-x_j^M$, respectively; a narrow ellipse $\partial E_S$ elongated along the $x_1$-axis
that intersects this axis at the points $(r_\pm^0,0)$; and a narrow ellipse $\partial E_L$ elongated along the $x_2$-axis
that intersects the $x_1$-axis at the points $(r_\pm,0)$. We denote by $E_S$ the closed domain bounded by the smaller ellipse, by $E_L$ the closed domain bounded by the larger ellipse,
and by $D_L$ and $D_S$ the discs bounded by the circles $\partial D_L$ and $\partial D_S$, see Fig.~\ref{Fig6} (left). That is, $E_S\subset D_S\subset D_L\subset E_L$. By construction,
the domain $E_S$ contains the critical points $(x_j^m,0), (x_j^M,0), (x_{j+1}^m,0)$ of the functional $V(\cdot,v^1)$, while all the other
critical points of this functional lie outside $E_L$. Let us denote by $J$ the segment $J=\{\bm{x}=(x_1,0): r_-\le x_1\le r_+\}$ and define the functional $\Phi$
by the formulas
\begin{equation}\label{new}
\Phi(x_1,x_2)=V(x_1,x_2,v^1) \quad {\rm for}\quad \bm{x}=(x_1,x_2)\in (\mathbb{R}^2\setminus E_L)\cup E_S \cup J,
\end{equation}
\begin{equation}\label{new1}
\Phi(x_1,x_2)=V\bigl(x_j^M+\sqrt{(x_1-x_j^M)^2+x_2^2},0,v^1\bigr)  \quad {\rm for}\quad \bm{x}=(x_1,x_2)\in D_L\setminus D_S
\end{equation}
(the latter formula is consistent with \eqref{sym}) in the union of the ellipse $E_S$, the annulus $D_L\setminus D_S$, the segment $J$
and outside the ellipse $E_L$. We extend this functional to the rest of the plane $(x_1,x_2)$  along vertical segments as follows.
For every vertical segment (of positive length) that belongs either to the upper half-plane $x_2\ge 0$ or the lower half-plane $x_2\le 0$
and that has one end $(x_1^e,x_2)$ on the circle $\partial D_S$ and the other end $(x_1^f,x_2)$
either on the curve $\partial E_S$ or the interval $\{(x_1,0): r_-^d\le x_1\le r_-^0\}$ or the interval $\{(x_1,0): r_+^0\le x_1\le r_+^d\}$,
the functional is defined by
\begin{equation}\label{phi3}
\Phi(x_1,x_2)=\frac{|x_2^2-(x_2^e)^2|}{|(x_2^e)^2-(x_2^f)^2|}\Phi(x_1,x_2^f)+\frac{|x_2^2-(x_2^f)^2|}{|(x_2^e)^2-(x_2^f)^2|}\Phi(x_1,x_2^e).
\end{equation}
The same formula defines $\Phi$ on  every vertical segment that belongs either to the upper half-plane $x_2\ge 0$ or the lower half-plane $x_2\le 0$
and that has one end $(x_1^e,x_2)$ on the curve $\partial E_L$ and the other end $(x_1^f,x_2)$
either on the circle $\partial D_L$ or the interval $\{(x_1,0): r_-\le x_1\le r_-^D\}$ or the interval $\{(x_1,0): r_+^D\le x_1\le r_+\}$.
By construction, the continuous and piecewise smooth functional $\Phi$ satisfies $\Phi(x_1,-x_2)=\Phi(x_1,x_2)$.
Furthermore, we assume that the ellipse $D_S$ has a sufficiently short vertical semi-axis $b_S$;
then relation \eqref{condi} implies that  $\Phi(x_1,x_2)$ is increasing in $x_2$
in the domain $\{(x_1,x_2): x_2\ge0\}\cap (D_S\setminus E_S)$.
We also assume that
the ellipse $D_L$ has a sufficiently
long vertical semi-axis $b_L$ so that $\Phi(x_1,x_2)$ is increasing in $x_2$  in the domain
$\{(x_1,x_2): x_1\in [r^D_-,r^D_+], x_2\ge0\}\cap (E_L\setminus D_L)$.
Since the functional $V(x_1,x_2,v^1)$ increases in $x_2$ in the upper half-plane
and \eqref{new} holds, we conclude that
$\Phi(x_1,x_2)$ increases in $x_2$ in the whole upper half-plane $x_2\ge0$ too.

Now we smoothen $\Phi$ using averaging over a small disc of radius $\rho$, that is we define
\begin{equation}\label{phit}
\tilde \Phi(x_1,x_2)=\frac{1}{\pi \rho^2} \iint\limits_{\xi_1^2+\xi_2^2\le \rho^2} \Phi(x_1+\xi_1,x_2+\xi_2)\,d\xi_1d\xi_2.
\end{equation}
This is a smooth functional with the derivatives defined by
\begin{equation}\label{derivatives}
\frac{\partial \tilde\Phi}{\partial x_1}(x_1,x_2)=\frac{1}{\pi \rho^2}\oint\limits_{\xi_1^2+\xi_2^2=\rho^2}
\Phi(x_1+\xi_1,x_2+\xi_2)\,d\xi_2,
\end{equation}
\begin{equation}\label{derivatives'}
\frac{\partial \tilde\Phi}{\partial x_2}(x_1,x_2)=- \frac{1}{\pi \rho^2}\oint\limits_{\xi_1^2+\xi_2^2=\rho^2}
\Phi(x_1+\xi_1,x_2+\xi_2)\,d\xi_1
\end{equation}
(with the integrals along the circle taken counterclockwise).
Since $\Phi$ is an even functional with respect to $x_2$, so is $\tilde \Phi$, that is
$\tilde \Phi(x_1,-x_2)=\tilde \Phi(x_1,x_2)$.
Therefore,
\begin{equation}\label{11}
\frac{\partial \tilde\Phi}{\partial x_2}=0\quad \ {\rm for}\quad \ x_2=0.
\end{equation}
Let us show that
\begin{equation}\label{12}
\frac{\partial \tilde\Phi}{\partial x_2}>0\quad \ {\rm for}\quad \ x_2>0.
\end{equation}
Indeed, if the circle of radius $\rho$ centered at the point $(x_1,x_2)$ belongs
to the upper half-plane $x_2\ge 0$, then relation \eqref{12} follows from
\eqref{derivatives'} because the functional $\Phi$ strictly increases in $x_2$ in this half-plane.
If part of the circle lies in the lower half-plane, then \eqref{12} follows from
the same equality due to the fact that the functional $\Phi$ strictly increases in $x_2$ and satisfies
$\Phi(x_1,-x_2)=\Phi(x_1,x_2)$.

Relations \eqref{11}, \eqref{12} imply that all the critical points of the functional $\tilde \Phi$ belong to the line
$x_2=0$ and are defined by the equation $\frac{\partial \tilde\Phi}{\partial x_1}=0$. Note that by construction $\tilde \Phi=\Phi$
in some neighborhood $\Omega$ of the critical points of the functional $\Phi$. Let us choose a sufficiently small parameter $\rho$ of smoothening in \eqref{phit} so that
the closed disc of radius $\rho$ centered at every critical point belongs to $\Omega$ and $\rho<r^D_+-r^d_+$.
Now, we prove that $\Phi$ and $\tilde \Phi$ have the same critical points. It suffices to show that
$\frac{\partial \tilde\Phi}{\partial x_1}(x_1,0)\ne 0$ on the intervals $[r_-,r^0_-]$ and $[r_+^0,r_+]$ (as on the rest of the $x_1$ axis
$\Phi=\tilde \Phi$).
To this end, let us take a $h\in[0,\rho]$ and consider the horizontal segment $J_h=\{(x_1,x_2): r_--\rho\le x_1\le r^0_-+\rho,\ x_2=h\}$.
We are going to show that $\Phi$ strictly decreases in $x_1$ along this segment and therefore equality \eqref{derivatives} (together with the relation $\Phi(x_1,-x_2)=\Phi(x_1,x_2)$)
implies that $\frac{\partial \tilde\Phi}{\partial x_1}(x_1,0)< 0$ for $x_1\in [r_-,r^0_-]$.

Consider the intersections $J_1$, $J_2$, $J_3$, $J_4$, $J_5$ of the interval $J_h$ with the sets
$\mathbb{R}^2\setminus E_L$, $E_L\setminus D_L$, $D_L\setminus D_S$, $D_S\setminus E_S$ and $E_S$, respectively (see Fig.~\ref{Fig6}).
On the intervals $J_1$ and $J_5$, the functional $\Phi$ decreases with the increasing $x_1$ because $\Phi=V(\cdot,v^1)=\tilde f(x_1)+h^2$
and $\tilde f'(x_1)<0$ in a neighborhood of the interval $[r_-,r_-^0]$.
On the interval $J_3$, according to \eqref{new1}, we have $\Phi(x_1,h)=V\bigl(x_j^M+\sqrt{(x_1-x_j^M)^2+h^2},0,v^1\bigr)=\tilde f(x_j^M+\sqrt{(x_1-x_j^M)^2+h^2})$,
which also implies $\frac{\partial \tilde\Phi}{\partial x_1}< 0$. It remains to consider the segments $J_2$ and $J_4$.
The segment $J_2$ can be divided into two parts,
$$
J_2^\ell=\bigl\{(x_1,h): x_j^M-a_L\sqrt{1-h^2/b_L^2}\le x_1 \le r_-^D\bigr\},
$$
$$
 J_2^r=\bigl\{(x_1,h): r_-^D \le x_1 \le x_j^M-\sqrt{R_L^2-h^2}\bigr\},
 $$
where $(x_1-x_j^M)^2/a_L^2+x_2^2/b_L^2=1$ is the equation of the ellipse $\partial E_L$ and
$(x_1-x_j^M)^2+x_2^2=R_L^2$ is the equation of the circle $\partial D_L$ (see Fig.~\ref{Fig6}).
According to \eqref{phi3}, on the interval $J_2^\ell$,
$$
\Phi(x_1,h)=\left(1-\frac{h^2}{b_L^2\bigl({1-{{a_L^{-2}}(x_1-x_j^M)^2}}\bigr)}\right)\tilde f(x_1)+\frac{h^2\bigl(\tilde f(x_1)+b_L^2\bigl(1-{a_L^{-2}}{(x_1-x_j^M)^2}\bigr)\bigr)}{b_L^2\bigl({1-{a_L^{-2}}{(x_1-x_j^M)^2}}\bigr)},
$$
hence $\Phi(x_1,h)=\tilde f(x_1)+h^2$ and
$\frac{\partial \Phi}{\partial x_1}(x_1,h)=\tilde f'(x_1)<0$ on $J_2^\ell$. On the interval $J_2^r$,
$$
\Phi(x_1,h)=\left(1-c(x_1)\right)\tilde f(x_j^M-R_L)+ c(x_1) \bigl(\tilde f(x_1)+ b_L^2({1-{a_L^{-2}}{(x_1-x_j^M)^2}})\bigr)
$$
with
$$
c(x_1)=\frac{h^2-(R_L^2-(x_1-x_j^M)^2)}{b_L^2(1-{a_L^{-2}}{(x_1-x_j^M)^2})-(R_L^2-(x_1-x_j^M)^2)},
$$
hence
$$
\begin{array}{rcl}
\frac{\partial \Phi}{\partial x_1}&=&c'(x_1)\bigl(\tilde f(x_1)+ b_L^2({1-{a_L^{-2}}{(x_1-x_j^M)^2}})-\tilde f(x_j^M-R_L)\bigr)\\
&+&c(x_1)\bigl(\tilde f'(x_1)-2(x_1-x_j^M)b_L^2 a_L^{-2}\bigr).
\end{array}
$$
Here
$$
\begin{array}{rcl}
{\rm sgn}\, c'(x_1)&=&{\rm sgn}\, \left[(b_L^2(1-{a_L^{-2}}{(x_1-x_j^M)^2})-(R_L^2-(x_1-x_j^M)^2))(x_1-x_j^M)\right.\\
&+&\left.(a_L^{-2} b_L^2-1)(x_1-x_j^M)(h^2-(R_L^2-(x_1-x_j^M)^2))\right]\ =\ -1
\end{array}
$$
because $b_L\gg 1$; and, $c(x_1)=O(h^2)$. Therefore, $\frac{\partial \Phi}{\partial x_1}<0$ for small $h$.
The segment $J_4$ can be divided into two parts,
$$
J_4^\ell=\bigl\{(x_1,h): x_j^M-\sqrt{R_S^2-h^2} \le x_1 \le r_-^0\bigr\},
$$
$$
J_4^r=\bigl\{(x_1,h): r_-^0\le x_1 \le x_j^M-a_S\sqrt{1-h^2 b_S^{-2}}\bigr\},
$$
where $(x_1-x_j^M)^2 a_S^{-2}+x_2^2 b_S^{-2}=1$ is the equation of the ellipse $\partial E_S$ and
$(x_1-x_j^M)^2+x_2^2=R_S^2$ is the equation of the circle $\partial D_S$. On the interval $J_4^\ell$,
$$
\Phi(x_1,h)=\left(1-\frac{h^2}{R_S^2-(x_1-x_j^M)^2}\right)\tilde f(x_1)+\frac{h^2}{R_S^2-(x_1-x_j^M)^2}\tilde f(x_j^M-R_S),
$$
and therefore,
$$
\frac{\partial \Phi}{\partial x_1}(x_1,h)=
\left(1-\frac{h^2}{R_S^2-(x_1-x_j^M)^2}\right)\tilde f'(x_1)+\frac{2h^2(x_1-x_j^M)}{(R_S^2-(x_1-x_j^M)^2)^2}(\tilde f(x_j^M-R_S)-\tilde f(x_1)),
$$
hence the relationships $\tilde f(x_j^M-R_S)>\tilde f(x_1)$ and $\tilde f'(x_1)<0$ imply
$\frac{\partial \Phi}{\partial x_1}<0$. Finally, on the segment $J_2^r$,
$$
\Phi(x_1,h)=\tilde c(x_1)\bigl(\tilde f(x_1)+b_S^2(1-{a_S^{-2}}{(x_1-x_j^M)^2})\bigr)+
(1-\tilde c(x_1))\tilde f(x_j^M-R_S)
$$
with
$$
\tilde c(x_1)=\frac{R^2_S-(x_1-x_j^M)^2-h^2}{R^2_S-(x_1-x_j^M)^2-b_S^2(1-{a_S^{-2}}{(x_1-x_j^M)^2})},
$$
hence
$$
\begin{array}{rcl}
\frac{\partial \Phi}{\partial x_1}(x_1,h)&=&
\tilde c'(x_1)\bigl(\tilde f(x_1)+ b_S^2({1-{a_S^{-2}}{(x_1-x_j^M)^2}})-\tilde f(x_j^M-R_S)\bigr)\\
&+&
\tilde c(x_1)(\tilde f'(x_1)-2(x_1-x_j^M)b_S^2 a_S^{-2}).
\end{array}
$$
Since $b_S\ll1$, $\tilde f'(x_1)<0$ and $f(x_j^M-R_S)>\tilde f(x_1)$, the term $\tilde c(x_1)(\tilde f'(x_1)-2(x_1-x_j^M)b_S^2 a_S^{-2})$ in this expression is negative
and the term $\bigl(\tilde f(x_1)+ b_S^2({1-{a_S^{-2}}{(x_1-x_j^M)^2}})-\tilde f(x_j^M-R_S)\bigr)$ is positive.
Furthermore,
$$
\begin{array}{rcl}
{\rm sgn}\, \tilde c'(x_1)&=&{\rm sgn}\, \bigl[-\bigl(R_S^2-(x_1-x_j^M)^2- b_S^2(1-{a_S^{-2}}(x_1-x_j^M)^2)\bigr)(x_1-x_j^M)
\\
&-&(R_L^2-(x_1-x_j^M)^2-h^2)(b_S^2 a_S^{-2}-1)(x_1-x_j^M)\bigr]\ =\ 1,
\end{array}
$$
where we again  use $b_S\ll 1$. Hence, $\frac{\partial \Phi}{\partial x_1}<0$ on this interval too.

We conclude that $\Phi$ decreases along each segment $J_h$, $|h|\le\rho$, and hence
$\frac{\partial \tilde\Phi}{\partial x_1}(x_1,0)< 0$ for $x_1\in [r_-,r^0_-]$. Similarly, one can show that
$\Phi$ strictly increases along each segment $\{(x_1,x_2): r_+^0-\rho\le x_1\le r_++\rho,\ x_2=h\}$ and therefore
$\frac{\partial \tilde\Phi}{\partial x_1}(x_1,0)>0$ for $x_1\in [r^0_+,r_+]$. Thus, the critical points of $\tilde \Phi$
coincide with the critical points of $\Phi$.

Now, we fix $v^2\in (v^1,u_+)$ and define the linear deformation
$$
V(\bm{x},u)=V(x_1,x_2,u)=(1-\tilde a(u)) V(x_1,x_2,v^1)+\tilde a(u)\tilde \Phi(x_1,x_2),\quad u\in[v^1,v^2],
$$
where $\tilde a(v^1)=0$, $\tilde a(v^2)=1$ and $\frac{d \tilde a}{du}(v^1)=\frac{d \tilde a}{du}(v^2)=0$.
Since both functionals $V(x_1,x_2,v^1)=\tilde f(x_1)+x_2^2$ and $\tilde \Phi(x_1,x_2)$ are even in $x_2$, strictly increase in $x_2$
for positive $x_2$, have the same critical points and the same monotonicity on the intervals of the line $x_2=0$ separated by the critical points,
this deformation also has the same critical points for every $u\in[v^1,v^2]$.
By construction (see formula \eqref{new1}), each circle $(x_1-x_j^M)^2+x_2^2=R^2$ within the annulus $D_L\setminus D_S$ is a level line of the functional $\Phi$
and since the smoothening parameter $\rho$ has been chosen small enough, there is an annulus ${\mathcal R}=\{(x_1,x_2): R_1< (x_1-x_j^M)^2+x_2^2< R_2^2\}\subset D_L\setminus D_S$ such
that each circle $(x_1-x_j^M)^2+x_2^2=R^2$ is a level line of the functional $\tilde \Phi$. Now, we can define a deformation
that interchanges the local minimum points $x_j^m=(x_j^m,0)$ and $x_{j+1}^m=(x_{j+1}^m,0)$ simply using the rotation:
\begin{equation}\label{rotate}
V(x_1,x_2,u)=\left\{
\begin{array}{ll}
V(x_1\cos \alpha -\sin\alpha, x_1\sin\alpha+x_2\cos\alpha,v^2),& (x_1-x_j^M)^2+x_2^2\le R,\\
V(x_1,x_2,v^2),& (x_1-x_j^M)^2+x_2^2> R
\end{array}
\right.
\end{equation}
with a smooth function $\alpha=\alpha(u)$ defined for $u\in[v^2,v^3]\subset [v^2,u_+)$, where $R\in (R_1,R_2)$ and $\alpha(v^2)=0$, $\alpha(v^3)=\pi$, $\frac{d\alpha}{du}(u)>0$ for $u\in(v^2,v^3)$
and $\frac{d\alpha}{du}(v^2)=\frac{d\alpha}{du}(v^3)=0$.
By construction, the functional $V(\cdot,v^3)$ increases in $x_2$ in the upper half-plane and is even in $x_2$. Finally, we use the linear deformation
$$
V(x_1,x_2,u)=(1-b(u)) V(x_1,x_2,v^3)+ b(u)(f(x_1)+x_2^2),\quad u\in[v^3,u_+],
$$
with $b(v^3)=0$, $b(u_+)=1$ to connect the functional $V(\cdot,v^3)$ to the functional $V(\cdot, u_+)=V(\cdot, u_-)$.
This completes the proof of the lemma and the theorem.

\section*{Appendix}
Here we briefly discuss multi-scale systems that may lead to equation \eqref{3} in the adiabatic limit.
Consider a two-dimensional system
$$
\ddot{\bm{x}}+ \gamma \dot{\bm{x}} +\nabla_x V(\bm{x},u)=0
$$
with $\bm{x}=(x_1,x_2)\in\mathbb{R}^2$ where the potential $V$ depends on a scalar parameter $u$. We assume large friction, $\gamma\gg 1$.
Introducing the slow time $\tau =\gamma^{-1} t$, this equation can be rewritten as a slow-fast system
$$
\begin{array}{rcl}
 \bm{x}' &=& \bm{y} \\
\varepsilon  \bm{y}' &=& - \nabla_x V(\bm{x},u)-\bm{y}
\end{array}
$$
with a small parameter $\varepsilon = \gamma^{-2}\ll 1$ (dot stands for $d/dt$, prime denotes the derivative  $d/d\tau$ with respect to the slow time). This system has an attracting slow manifold which approaches the critical surface $\bm{y} = - \nabla_x V(\bm{x},u)$ in the limit $\varepsilon\to 0$
(equivalently, $\gamma\to \infty)$. Hence, for large $\gamma$, dynamics on this manifold can be approximated by the gradient equation
\begin{equation}\label{3'}
\bm{ x}' = -\nabla_x V(\bm{x},u)
\end{equation}
or, $\gamma \dot{\bm{x}} = -\nabla_x V(\bm{x},u)$ using the original time $t$.

Now, we assume that the parameter $u$ can vary in time slowly, that is $u=u(\nu t)$.
Hence, we have a system with three time scales,
\begin{equation}\label{4}
\ddot{\bm{x}}+ \gamma \dot{\bm{x}} +\nabla_x V(\bm{x},u(\nu t))=0,
\end{equation}
and the approximating system with two time scales:
\begin{equation}\label{5}
\bm{x}' = -\nabla_x V(\bm{x},u (\mu \tau))
\end{equation}
where $\mu=\gamma \nu$.
The variable $u$ can be interpreted as input of these systems.
Finally, assume that the characteristic time scale associated with the input variations, $\theta=\nu t$,
is the slowest time scale and, furthermore, $\mu=\gamma \nu\ll 1$.
In this case, dynamics of system \eqref{5} (as well as dynamics of system \eqref{4} in the limit
$\gamma\to \infty$, $\mu=\gamma \nu\to 0$) on the slowest time scale of the input variations, $\theta=\nu t$, can be described as a sequence of transitions
between local minimum points of the  potential energy $V(\cdot,u)$ along the gradient flow of system \eqref{3'}.
%
We remark that hysteresis is a cause of energy dissipation. Namely, each 
transition from
a critical point $x_*(u)$ to a critical point $x^*(u)$ of the potential $V(\cdot,u)$ (induced by a saddle-node bifurcation that eliminates the critical point $x_*(u)$)
is associated with 
irreversible dissipation of  $V(x^*(u),u)-V(x_*(u),u)$ units of energy.
Continuous evolution of the system with a minimum point $x_*(u)$ of the potential
is reversible and dissipationless.


%
%


\end{document}